\documentclass[12pt]{amsart}
\usepackage[T1]{fontenc}

\title{Les nombres de Cuesta-Conway comme extension des ordinaux de Cantor\\ \tiny Une courte introduction aux nombres surr\'eels} 

\author{Labib Haddad}
\address{120 rue de Charonne, 75011 Paris, France}
\email{labib.haddad@wanadoo.fr}

\usepackage{amssymb}
\usepackage{amsmath}
\usepackage{endnotes}

\newcommand{\su}{\subsection*}
\newcommand{\head}{\section*}
\newcommand{\noi}{\noindent}

\newcommand{\Ž}{\'e}

\newcommand{\ˆ}{\`a}
\newcommand{\}{\`u}

\newcommand{\Q}{\mathbb Q}
\newcommand{\R}{\mathbb R}

\newcommand{\sO}{\mathsf O}

\newcommand{\sS}{\mathsf S}

\newcommand{\leqs}{\leqslant}

\newcommand{\geqs}{\geqslant}

\newcommand{\guil}{\guillemotleft}  
\newcommand {\guir}{\guillemotright}

\newcommand{\ali} {\begin{aligned}}   
\newcommand{\ala} {\end{aligned}}

\newcommand {\et}{\ \text{et}\ }

\newcommand {\ou}{\ \text{ou}\ }

\newcommand {\pour}{\ \text{pour}\ }

\newcommand{\stm}{\smallsetminus}

\newcommand{\vvide}{\varnothing}

\newcommand{\inc}{\subset}

\newcommand{\bc}{\begin{cases}}
\newcommand{\ec}{\end{cases}}

\newcommand{\lopar}{\noi \{$\looparrowright$ \ }

\begin{document}
\maketitle

\thispagestyle{empty}

\markboth{Labib Haddad}{ Nombres de Cuesta-Conway}

\

\head{Pr\Žambule}

\

\

Avec 30 ans d'\Žcart, en 1954 puis 1976, par deux proc\Žd\Žs diff\Žrents, ind\Žpendamment l'un de l'autre, {\sc Cuesta} [ 2 ]  et {\sc Conway}  [ 3 ] ont introduit la m\me classe de nombres,  extension de la classe  des ordinaux de {\sc Cantor}.

\

On d\Žsignera par $\sS$ cette classe de nombres de {\sc Cuesta-Conway} et par $\sO$ la classe des nombres ordinaux de {\sc Cantor}. Pour abr\Žger, tout \Žl\Žment de la classe $\sS$ sera appel\Ž {\bf nombre}, sans autre qualificatif, et tout \Žl\Žment de la classe $\sO$  sera appel\Ž  {\bf  ordinal}. 

\

{\sc D. Knuth} avait dit, en son temps, que les nombres de Conway devraient s'appeler {\bf les nombres surr\Žels}. C'est pourquoi on a choisi de d\Žsigner  cette classe de nombres par le symbole $\sS$, tout nu, sans fioriture !

\

\lopar Bien entendu, c'est  {\sc Conway} qui a d\Žvelopp\Ž au del\ˆ de la simple d\Žfinition la th\Žorie de cette classe de nombres qu'il d\Žsigne {\bf $\mathbf {No}$}. Cependant, il n'est peut-\tre pas inutile d'associer le nom de {\sc Cuesta} pour lui rendre ainsi hommage.  Cela n'enl\ve rien \ˆ l'immense m\Žrite de {\sc Conway} !\}

\

Il y a plusieurs mani\res de d\Žfinir cette classe de nombres de {\sc Cuesta-Conway}. Certaines utilisent la notion de coupure de Dedekind dans le corps $\Q$ des nombres rationnels, \Žtendue aux ensembles totalement ordonn\Žs quelconques.  

\

J'ai d\Žj\ˆ eu l'occasion de parler des constructions de {\sc Cuesta} et de {\sc Conway} dans un article post\Žs sur {\tt arXiv}, voir [ 5 ]. Je ne reviens pas l\ˆ-dessus.

\

Outre ces constructions, on  en  trouve une dans  {\sc Gonshor} [ 4 ], par exemple, et plusieurs  autres signal\Žes dans {\sc Alling} [ 1 ].

\

Je voudrais pr\Žsenter, ici, une nouvelle mani\re, assez directe, pour d\Žfinir et construire la classe $\sS$ des nombres,  qui \Žvite de longs d\Žtours avant d'arriver \ˆ {\bf la forme normale}, en allant droit au but ! J'utilise pour cela la notion de sommes formelles index\Žes par les ordinaux. Voici de quoi il s'agit !

\

\

\head{La classe des nombres}
\

La classe $\sO$ est une sous-classe de la classe $\sS$. Autrement dit, les ordinaux sont des nombres : $\sO \inc \sS$. La classe $\sS$ est munie d'une addition, \guil $ \ + \ $\guir, et d'une multiplication, \guil $  \ \cdot \ $\guir , qui co\•ncident avec l'addition et la multiplication des ordinaux, pour les ordinaux.

\

\

{\bf La classe des nombres est un corps commutatif totalement ordonn\Ž qui renferme les ordinaux et le corps des nombres r\Žels $\R$.}

\

La classe des nombres est totalement ordonn\Že \ˆ l'aide d'une relation d'orde d\Žsign\Že par $\leqs$. On d\Žsigne par  $<$ la relation d'ordre stricte associ\Že \ˆ $\leqs$. Bien entendu, $\geqs$ et $>$, d\Žsignent les oppos\Žes.

\

 Dans $\sS$,  \ˆ chaque nombre $y$ est associ\Ž un nombre  $\omega^y$, appel\Ž {\bf son exponentielle}. Lorsque l'on a $\alpha\in \sO$, l'exponentielle, $\omega^\alpha$, co\•ncide avec l'ordinal $\omega^\alpha$.

\

L'exponentielle poss\de les propri\Žt\Žs suivantes :
$$\omega^0 = 1$$
$$ \omega ^ {a + b} = \omega ^ a. \omega ^ b$$
$$a < b \implies \omega^a < \omega^b.$$

\

Les nombres sont les sommes formelles, index\Žes par des ordinaux, de la forme suivante :
\[a = \sum_{i < \alpha} \omega^{a_i} r_i,\tag*{(1)}\]
o\ $\alpha$ est un ordinal  donn\Ž quelconque, $(a_i)_{i < \alpha}$ une suite strictement d\Žcroissante de nombres, et $(r_i)_{i < \alpha}$ une suite de nombres r\Žels non nuls. Les deux suites sont, toutes deux,  index\Žes par les ordinaux $i < \alpha$.

\

\lopar C'est ce que l'on appelle la {\bf forme normale} du nombre $a$. Elle co\•ncide avec la forme normale de Cantor pour les ordinaux. Voir {\sc Alling} [ 1 ], p.258, section  7.11.\}

\

Pour $\alpha = 0$, la somme n'a aucun terme;  elle est \Žgale \ˆ l'entier $0$, par convention.

\

Pour faire court, on dira que le nombre $a$, dans la formule (1), est de {\bf hauteur} $\alpha$ et on \Žcrira $h(a) = \alpha$.

\

L'ordre total sur $\sS$  est d\Žfini comme suit. Il doit \tre compatible avec l'addition et, pour 
$$0\neq a = \sum_{i< \alpha} \omega^{a_i} r_i,$$
on a
$$a > 0 \iff r_0>0.$$

\

\noi Cet ordre total  est ainsi lexicographique. Plus pr\Žcis\Žment, soient
$$a = \sum_{i<\alpha} \omega^{a_i} r_i \neq b = \sum_{i < \beta} \omega^{b_i} s_i,$$
deux nombres non nuls et
soit $k$ le plus petit des $i$ pour lesquels on a
 $$\omega^{a_i} r_i \neq  \omega^{b_i} s_i.$$
 Alors, on aura $a < b$ 
si ($a_k = b_k$, et $r_k < s_k$) ou ($a_k < b_k$).

\

\noi Les ordinaux correspondent aux cas o\, pour tout $i$ dans (1), le $a_i$ est un ordinal et le $r_i $  est un entier naturel.  Les nombres r\Žels correspondent aux cas o\  ($\alpha = 1$ et $a_0 = 0$), et au cas particulier $\alpha = 0$, bien entendu.

\

\su{Reformulation} Tout cela \Žtant dit, on observe que l'essentiel dans la d\Žfinition de  ces sommes formelles n'est pas l'usage du segment initial d'ordinaux, $[0,\alpha[$ . On pourrait  le remplacer  par un ensemble bien ordonn\Ž quelconque !  On exige seulement que la suite strictement d\Žcroissante des $a_i$ ait pour type l'ordre  oppos\Ž \ˆ celui d'un segment initial  d'ordinaux, $[0,\alpha[$.

\

On reformule alors, afin d'abr\Žger et de simplifier les \Žcritures. On appellera {\bf couple ad\Žquat} tout couple $(A,r)$ o\ $A$ est un sous-ensemble (vide ou non) de la classe des nombres, ($A \inc \sS$), tel que $(A, >)$ soit bien ordonn\Ž  et o\ $r = (r_x)_{x\in A}$ est  une famille de nombres r\Žels non nuls. Il faut insister sur les deux mots, {\bf sous-ensemble} et {\bf non nuls} !

\

Ainsi, \ˆ chaque couple ad\Žquat est associ\Že une forme normale et vice  versa !

\

On pourra alors \Žcrire le m\me nombre, $a$, sous deux formes diff\Žrentes :

\

$a = (A,r)$, sa {\bf forme courte}, comme couple ad\Žquat,

\

$a = \sum_{i<\alpha} \omega^{a_i} r_i$, sa {\bf forme normale}, ayant pour hauteur $h(a)=\alpha$.

\

\lopar De mani\re analogue, on \Žcrit bien un nombre rationnel sous plusieurs  formes diff\Žrentes : par exemple, comme fraction $p/q$, ou comme d\Žveloppement d\Žcimal  p\Žriodique \ˆ partir d'un certain rang.\}

\

\head{Les formes longues}

\

Il sera parfois commode d'\Žcrire le nombre $a$ sous {\bf une forme longue}
$$\sum_{i <\gamma}\omega^{a_i} r_i$$
o\ $(a_i)_{i<\gamma}$ est une suite  d\Žcroissante de nombres, pas n\Žcessaiement strictement d\Žcroissante,  et o\ $r= (r_i)_{i<\gamma}$ est une suite de nombres r\Žels dont certains peuvent \tre nuls. La suite    d\Žcroissante, $(a_i)_{i<\gamma}$, peut marquer des paliers, avec des termes successifs \Žgaux : 

$a_i= a_{i +1} = \dots = a_{i + n} = \cdots$.

\noi On exige cependant, chaque fois, que ce  palier ait une longueur finie. Comme pour les polyn\™mes, on peut alors regrouper   les termes
$$\omega^{a_i}  r_i + \dots + \omega^{a_{i + n}}r_{i + n} = \omega^{a_i}(r_i + \dots + r_{i + n})$$
puis  \Žliminer, le cas \Žch\Žant, les termes dont les coefficients r\Žels sont nuls ! On obtient une somme que l'on peut r\Žindexer \ˆ l'aide des ordinaux d'un segment initial $[0,\alpha[$ et obtenir  un nombre $a$ sous forme normale.

\

\lopar On notera que, pour un nombre donn\Ž, sa forme normale est la moins longue de toutes ses formes longues !  

Parmi  les formes longue, on distingue une classe particuli\re, celles dont les suites $(a_i)$ sont {\bf strictement} d\Žcroissantes.\}

\

\noi On dispose, de cette fa\c con,  de plusieurs mani\res pour \Žcrire le m\me nombre :
forme normale, forme courte, formes longues. On utilisera l'une ou l'autre de ces formes, selon le cas. On essaiera de choisir, chaque fois,  la plus commode !

\

\lopar Toutes ces sommes formelles, formes normales, formes longues, sont  connues sous le nom de {\bf s\Žries de Hahn}. Voir l'article [ 7 ],  par exemple, pour la d\Žfinition.

\

\noi Les nombres se pr\Žsentent ainsi comme des polyn\™mes g\Žn\Žralis\Žs. On  le verra ci-dessous, pour l'addition et la multiplication, ils se comportent comme des polyn\™mes ou des s\Žries formelles classiques. Ainsi, la classe $\sS$  est une sorte de $\R$-alg\bre g\Žn\Žralis\Že, engendr\Že par des exponentielles !

\

\noi On utilisera, ci-dessous, les formes longues en particulier lorsque l'on voudra que des sommes diff\Žrentes soient index\Žes par les m\mes ordinaux et, le cas \Žch\Žant,  qu'elles aient la m\me suite $ (a_i)_{i < \gamma}$ ! Ce proc\Žd\Ž permettra de simplifier beaucoup  de d\Žfinitions et de rendre les calculs plus ais\Žs \ˆ suivre.\}

\

\head{L'addition}  

\

Soient 
$$a = \sum_{i< \gamma} \omega^{a_i}r_i  \ , \ b = \sum_{i <\gamma} \omega^{a_i}s_i$$ 
deux nombres sous formes longues. On pose
$$a + b = \sum_{i<\gamma} \omega^{a_i} (r_i + s_i).$$
Cela d\Žfinit l'addition pour les nombres, sous formes longues.

\

\lopar Sans grand d\Žtour, on peut s'assurer que cette d\Žfinition ne d\Žpend  que des nombres $a$ et $b$, pas de leurs formes longues choisies !\}

\

On peut v\Žrifier simplement que la classe $\sS$ munie de l'addition et de la relation d'ordre totale ainsi d\Žfines est un groupe commutatif, totalement ordonn\Ž ayant $0$ pour \Žl\Žment  neutre.

\

\lopar {\bf  Rappels sur les notations et le vocabulaire des relations d'ordre}

$a\geqs b$ se lit : $a$ est plus grand que  $b$.

$a >  b$ se lit : $a$ est strictement plus grand que $b$.

$a \geqs 0$ se lit : $a$ est positif.

$a > 0$ se lit : $a$ est strictement positif.\}

\

On se servira d'un lemme de {\sc H. B. Neumann} [ 8, p.204-209 ] dont on trouve une d\Žmonstration relativement courte dans [ 6 ]. Il s'agit d'un lemme combinatoire dont on utilisera le cas particulier relatif aux groupes commutatifs totalement ordonn\Žs.

\

\head{Le lemme de H. B. Neumann}

\

Voici la version particuli\re du lemme de {\sc Neumann} dont nous nous servirons.

\

{ \it Dans le groupe commutatif totalement ordonn\Ž $\sS$, soit $A$ une partie  bien ordonn\Že pour la relation $\geqs$, form\Že d'\Žl\Žments strictement positifs. L'ensemble des sommes d'un nombre fini d'\Žl\Žments de $A$ est lui-m\me bien ordonn\Ž pour la relation $\geqs$. De plus, pour chaque nombre $t$, il n'y a qu'un nombre fini de mani\res d'\Žcrire $t$ comme somme d'un nombre fini d'\Žl\Žments de $A$.}

\

Dans la suite, on se servira souvent (au moins quatre fois) de ce lemme.

\su{Convention}  \`A chaque ensemble $A$ de nombres strictement positifs tel que $(A,\geqs)$ soit bien ordonn\Ž, on associe l'ensemble des sommes d'un nombre fini d'\Žl\Žments de $A$ que l'on appellera {\bf la doublure} de l'ensemble$A$.
 Ordonn\Že par  $\geqs$, cette doublure est bien ordonn\Že, d'apr\s le lemme de Neumann.

\

\head{La multiplication}

\

On se donne deux nombres sous formes longues
$$a = \sum_{i<  \alpha} \omega^{a_i} r_i  \ , \ b = \sum_{j < \beta} \omega^{b_j} s_j.$$ 
On leur associe la famille de {\bf  mon\™mes} suivante : 
$$\left\{\omega^{a_i + b_j}r_i. s_j\right\}$$
analogues aux mon\™mes dans la multiplication de deux polyn\™mes ! 
On introduit l'ensemble $C = \{a_i + b_j : i < \alpha,\ , \ j < \beta\}$.
C'est une partie de l'ensemble des sommes d'un nombre fini d'\Žl\Žments de l'ensemble bien ordonn\Ž, $(A\cup B, \geqs)$, autrement dit, c'est une partie de la doublure de l'ensemble $A\cup B$, de sorte que $(C, \geqs)$ est bien ordonn\Ž. 

 \
 
 \noi On indexe l'ensemble $C$ en une suite strictement d\Žcroissante $(c_k)_{k< \gamma}$. 
 
 \
 
 \noi De plus, pour chaque $k$ donn\Ž, il n'y a qu'un nombre fini de mani\res d'\Žcrire  $c_k = a_i + b_j$. Soit
 $F(k) = \{(i, j) : c_k = a_i +  b_j , \ i < \alpha \ , \ j <\beta\}$. On pose
$$t_k = \sum_{(i,j) \in F(k)} r_i.s_j.$$
C'est une somme d'un nombre fini de nombres r\Žels, c'est donc un nombre r\Žel qui peut \tre nul, le cas \Žch\Žant ! La somme formelle  
$$c = \sum_{k < \gamma} \omega^{c_k} t_k.$$
est la forme longue d'un nombre $c$ dont la suite des $c_k$ est strictement d\Žcroissante.

\

\noi On remarquera que le nombre $c$ qui d\Žpend des deux nombres, $a$ et $b$, ne d\Žpend pas du choix de leurs formes longues puisque  l'on a $r_i.s_j = 0$ si $r_i = 0 \ou s_j = 0$.

\

Le nombre 
$$c = \sum_{k < \gamma} \omega^{c_k} t_k$$
construit ci-dessus \ˆ partir des deux nombres donn\Žs, $a$ et $b$, est le produit de ces deux nombres, {\bf par d\Žfinition} ! C'est ainsi que l'on d\Žfinit la multiplication des nombres. Autrement dit, on a
$$a.b = c = \sum_{k < \gamma} \omega^{c_k} t_k.$$

 \

\head{V\Žrifications}

\

On v\Žrifie que la classe $\sS$ munie de l'addition, de la multiplication et de la relation d'ordre total ainsi d\Žfinies, est un anneau commutatif, totalement ordonn\Ž, ayant $1$ pour \Žl\Žment unit\Ž. Les v\Žrifications se font sans grand d\Žtour, mais elles sont {\it longues} et {\it fastidieuses} !

\

\head{L'anneau des nombres est un anneau int\gre} 

\

Soient $a = (A,r) \neq 0 \ , \ b = (B,s) \neq 0 \ , \ c = a.b = (C,t)$ sous formes courtes.
Dire que le nombre $a= (A,r)$ n'est pas nul c'est dire que $A$ n'est pas vide. Donc $A$ poss\de un \Žl\Žment maximal : on le d\Žsignera par $a_0$. De m\me, le nombre $b= (B,s)$ n'\Žtant pas nul, l'ensemble $B$ poss\de un \Žl\Žment maximal  $b_0$. Lorsqu'il en est ainsi, l'ensemble $C$ n'est pas vide et son \Žl\Žment maximal est $c_0 = a_0 + b_0$ et l'on aura
$t_{c_0} = r_{a_0}.s_{b_0}$.
Ainsi
$$a \neq 0 \et b \neq 0 \iff a.b \neq 0.$$
Autrement dit, l'anneau $\sS$ est un anneau int\gre !

\

Pour montrer que $\sS$ est un corps, il reste \ˆ prouver que tout nombre non nul poss\de un inverse

\

\head{Chaque nombre non nul poss\de un inverse}

\

On montre que, pour un nombre donn\Ž quelconque $a \neq 0$, il existe toujours un nombre $b$ tel que $a.b = 1$. Autrement dit, il existe un nombre $b$ inverse
du nombre $a$, ce que l'on \Žcrit $b = a^{-1} = 1/a$. Cet inverse est bien entendu unique puisque $\sS$ est un anneau int\gre !

\

\su{D\Žmonstration} Tout nombre $a = \sum_{i < \alpha}\omega^{a_i} r_i$ non nul se met sous la forme :
$$a = \omega^{a_0} r_0\left(1 +\sum_{i< \alpha} \omega^{a_i-a_0} r_i/r_0\right)$$
et $\omega^{-a_0}/ r_0$ est l'inverse de   $\omega^{a_0} r_0$. Il suffit donc de d\Žmontrer le r\Žsultat pour les nombres dont la forme normale commence par $1$. On proc\de d'une mani\re analogue \ˆ celle que l'on emploie pour les s\Žries enti\res classiques \ˆ coeficients r\Žels : par approximations successives et r\Žcurrence.

\

\noi On se donne le nombre $a$ sous forme normale : 
$$a = 1 + \sum_{i <\alpha} \omega^{a_i} r_i.$$
On sait que $(A =\{a_i : i<\alpha\}, \geqs)$ est bien ordonn\Ž. On d\Žsigne par $B= \{b_j : j < \beta\}$ la doublure de $A$. On montre qu'il y a un nombre $b$ sous forme longue : 
$$b = 1 +  \sum_{j < \beta} \omega^{b_j} s_j, \ \text{o\} \  \{b_j : j < \beta\} = B,$$  
tel que $a.b = 1$.  Pour cela, on \Žcrit
$$a.b = c = 1 + \sum_{k < \gamma} \omega^{c_k} t_k$$
sous sa forme longue o\ la suite $(c_k)$ est strictement d\Žcroissante, et le coefficient r\Žel $t_k$ est donn\Ž par la formule suivante
$$t_k = \sum_{(i,j) \in F(k)} r_is_j \ , \ F(k) = \{(i,j) : c_k = a_i + b_j\}\ \text{qui est fini}.$$
Puisque les suites $(a_i), (b_j), (c_k)$ sont strictement d\Žcroissantes, on aura
$$a_i  < 0 \ , \ b_j < 0  \ , \ c_k < 0, \pour i >0, j>0, k > 0.$$
Dans le cas pr\Žsent, on a $a_0 = b_0 = c_0 = 0$ et $r_0 = s_0 = t_0 = 1$. 

\noi On pose $(b_1,s_1) = (a_1, -r_1)$ et on prend $b = 1+\omega^{b_1} s_1$,  on aura ainsi $c_1 = a_1 + 0 = 0 + a_1$ donc  $t_1 = r_1 + s_1= 0 $ d'o\ $t_1 = 0$.

\

\noi On suppose que, pour  $b = 1 + \sum_{j < k} \omega^{b_j} s_j$, on a $t_m = 0$, pour $0< m < k$. Alors, en prenant  $s_k = -\sum_{(i,j) \in F(k)\stm (0,k)} r_is_j$ et $b = 1 + \sum_{j \leqs k} \omega^{b_j} s_j$, on aura $t_m = 0, \pour 0 < m \leqs k$. Ainsi de suite pour tout $k < \beta$.  Avec le nombre $b = 1 + \sum_{j < \beta} \omega^{b_j} s_j$ ainsi construit, on aura $ab=c = 1$.\qed

\

\

{\bf On a bien montr\Ž que la classe des nombres est un corps commutatif totalement ordonn\Ž qui renferme les ordinaux et le corps des nombres r\Žels $\R$.}

\

\

\head{Construction de la classe des nombres}

\

On  construit la classe des nombres \guil \ couche  apr\s couche \guir, par r\Žcurrence transfinie. On d\Žfinit une suite croissante d'ensembles, $(K_\alpha)_{\alpha \in \sO}$  dont la r\Žunion est  $\sS$, la classe des nombres. 

\

\`A l'\Žtape $\alpha$ de la r\Žcurrence, on pose 

$E_\alpha = \bigcup_{\beta < \alpha}K_\beta,$

$K_\alpha = \text{ensemble des nombres $a = (A,r)$ tels que $A\inc E_\alpha$}$.

 On part de $K_{-1} = \vvide$  de sorte que $K_0  = \{0\}$.

\

La r\Žcurrence est faite ! On a 
$$\sS = \bigcup_{\alpha \in \sO} K_\alpha.$$

\

\noi En se souvenant de  ce que sont l'addition, la multiplication et  l'inverse d'un nombre, il est clair  que chacun des $K_\alpha$ ($\alpha > 0$) est un corps ! Ce sont des sous-corps du corps des nombres $\sS$. En particulier, $K_1$ est le corps des nombres r\Žels ! Pour les nommer, on dira que ces $K_\alpha$ sont {\bf  les corps liminaires}

\

\lopar  Il est bon de redire que $\sS$ est  {\bf une  classe}, pas {\bf un ensemble} !\}

\

\head{Le corps des nombres, $\sS$, est r\Žellement clos}

\

Cela veut dire qu'aucune de ses extensions alg\briques propres n'est totalement ordonnable. C'est une des propri\Žt\Žs remarquables du  corps des nombres, $\sS$. D'apr\s un th\Žor\me classique, un corps commutatif totalement ordonn\Ž est r\Žellement clos si et seulement s'il poss\de les deux propri\Žt\Ž suivantes.

\su{1} Chaque \Žl\Žment positif y poss\de  une racine carr\Že.

\su{2} Tout polyn\™me de degr\Ž impair \ˆ coefficients dans le corps y poss\de au moins une racine.

\

\noi On va montrer, plus g\Žn\Žralement, que  tous les corps liminaires, $K_\alpha$,  poss\dent  ces deux propri\Žt\Žs. Ils sont donc, tous, des corps r\Žellement clos. Il en d\Žcoule que le corps des nombres, $\sS$, lui-m\me, r\Žunion croissantes des corps liminaires,  est r\Žellement clos. Le premier corps liminaire, $K_1$, n'est autre que le corps $\R$ des nombres r\Žels qui, on le sait, est r\Žellement clos.

\

\su{Proposition} {\sl Dans un coprs liminaire, $K_\alpha$, chaque nombre positif poss\de une racine carr\Že.} 

\su{D\Žmonstration} Soit $a = \sum_{i < \alpha}\omega^{a_i}r_i$ un nombre strictement positif, sous forme normale. On peut \Žcrire
$$a = \omega^{a_0}r_0 \left(1 + \sum_{i < \alpha}\omega^{a_i-a_0}\frac{r_i}{r_0}\right).$$

Le nombre $\omega^{a_0}r_0$ poss\de la racine carr\Že $\omega^{a_0/2}\sqrt {r_0}$.
Il suffit donc de consid\Žrer les formes normales qui commencent par $1$. 

On se donne $a = 1 + \sum_{i < \alpha}\omega^{a_i}r_i$ sous forme normales. On d\Žsigne par $B =\{b_j : j < \beta\}$ la doublure de $A = \{a_i : i<\alpha\}$.  On pose $b = 1 + \sum_{j < \beta} \omega^{b_j}s_j$. On d\Žtermine les inconnues $s_j$, par r\Žcurrence, afin d'avoir $b^2 = a$. 

\

\noi On suppose que $c$ est le premier exposant pour lequel les coefficients de $\omega^c$ diff\rent dans $a = \sum_{i < \alpha}\omega^{a_i}r_i$ et dans $(1 + \sum_{i < \beta} \omega^{b_i}s_i)^2$. Soient $r$ et $t$, respectivement, ces coefficients de $\omega^c$ dans $a$ et  dans $b^2$. Si l'on ajoute \ˆ $b$ le terme $\omega^c d$, le coefficient de $\omega^c$ dans $b^2$ devient $t + 2d$, bien entendu. En prenant $2d = r-t$, le coefficient de $\omega^c$ dans $b^2$ devient \Žgal \ˆ $r$. En poursuivant la r\Žcurrence, le nombre $b$ que l'on obtient ainsi est tel que $b^2 =a$.\qed 

\

\su{Th\or\me} {\sl Dans tout corps liminaire, $K_\alpha$, chaque polyn\™me de degr\Ž impair poss\de au moins une racine.} Bien entendu, il suffit d'\Žtablir le r\Žsultat pour les polyn\™mes irr\Žductibles car un polyn\™me de degr\Ž impair poss\de toujours un  facteur irr\Žductible de degr\Ž impair  !

\

\lopar  Dans [ 4 , p.73-76, lemme 5.7 et th\Žor\me 5.10 ], {\sc  Gonshor} \Žtablit ce r\Žsultat pour le corps des nombres, $\sS$.  N'ayant pas fait mention des corps liminaires,  $K_\alpha$, il n'en parle pas, bien entendu ! La d\Žmonstration que l'on donne ici s'inspire {\bf largement} de celle de {\sc Gonshor} ! Avec quelques petites retouches !\}

\

\head{D\Žmonstration du th\Žor\me}

\

Elle comporte une d\Žfinition, un lemme et une fin.

\su{D\Žfinition} Un {\bf infinit\Žsimal}  est un nombre  $t$ tel que $-r < t < r$ pour tout nombre r\Žel $r > 0$.

\
 
\head{Un lemme \emph{\ˆ la Hensel}}

\

On \Žtablit le r\Žsultat de {\it rel\vement  \ˆ la Hensel} suivant.

\

Soit $f(x) = x^n + \sum_{i = 1}^n c_ix^{n-i}$ un polyn\™me dont les coefficients sont des nombres de la forme $c_i = r_i + t_i$ o\ $r_i$ est un nombre r\Žel et $t_i$ un infinit\Žsimal. On suppose que le polyn\™me $P(x) = x^n + \sum_{i=1}^n r_ix^{n-i}$ \ˆ coefficients r\Žels, {\bf la partie r\Želle} du polyn\™me $f(x)$, est le produit de deux polyn\™mes, $Q$ et $R$, premiers entre eux, ayant pour degr\Žs $q$ et $r$, respectivement. Alors $f$ est le produit de deux polyn\™mes $g$ et $h$ ayant pour degr\Žs, $q$ et $r$, respectivement.
\su{D\Žmonstration} En \Žcrivant les infinit\Žsimaux, $t_i$, sous leurs formes longues, on met le polyn\™me $f(x)$ sous la forme
$$f(x) = \sum_{i < \alpha} \omega^{a_i} s_i$$
o\ $s_i$ est  un polyn\™me \ˆ coeffcients r\Žels de degr\Ž au plus \Žgal \ˆ $n-1$, pour $i > 0$,  $s_0 = x^n + \sum_{i=1}^ n r_ix^{n-i}$ et $a_0 = 0$.

\

\noi Ici, tous les $a_i$ sont strictement n\Žgatifs puisque les $t_i$ sont infinit\Žsimaux. Par la sym\Žtrie $x \mapsto -x$, on obtient $A = \{-a_i : i <   \alpha\}$ bien ordonn\Ž, et sa doublure, $B = \{- b_i : i < \beta\}$. Comme pour la racine carr\Že, on intoduit
$$U  = \sum_{i< \beta} \omega^{b_i}Q_i \ , V=\sum_{i< \beta}\omega^{b_i}R_i.$$
Puis on d\Žtermine les suites de  polyn\™mes, $Q_i, R_i$, afin que l'on ait
 $$UV = \left(\sum_{i< \beta} \omega^{b_i}Q_i\right)\left(\sum_{i< \beta}\omega^{b_i}R_i\right) =f(x)$$
 avec $Q_0 = Q, R_0 = R$  et, pour  $i > 0$, $\deg(Q_i) < q$, $\deg(R_i) < r$.
 
 \
 
 \noi Soit $c$ le premier exposant pour lequel les coefficients de $\omega^c$ diff\rent dans $f(x)$ et dans  le produit $UV$. Soient $W$ et $T$, respectivement, ces deux coefficients. Le coefficient de $\omega^{a_0}$ dans $f(x)$ et dans le produit $UV$ est le m\me, \Žgal \ˆ $QR$, de sorte que l'on a $c > a_0$. [De  fait, $W$ est le polyn\™me $s_i$, pour $a_i = c$.] On ajoute \ˆ $U$ et \ˆ $V$, respectivement, les termes  $\omega^c G$  et $\omega^c H$,  o\  $G$  et $H$ sont des polyn\™mes \ˆ d\Žterminer.  On obtient
 $$(U + \omega^c G)(V + \omega^c H) = UV + \omega^c(HQ + GR) + \cdots.$$
 Le coefficient de $\omega^c$ devient  $T + (HQ  + GR)$ dans ce produit.  
On veut que ce coefficient soit \Žgal \ˆ celui de $\omega^c$ dans $f(x)$. Il suffit donc de trouver des polyn\™mes, $G$  et $H$, de sorte que l'on ait
 $$HQ+ GR= W-T \et \deg(G) \leqs q-1, \deg(H)\leqs r-1.$$
Or, $Q$ et $R$  satisfont une identit\Ž de B\Žzout, $XQ + YR = 1$, puisqu'ils sont premiers entre eux, de sorte que l'on a
$$(W-T)XQ + (W-T)YR   = W-T.$$
En posant 
$$E = (W-T)X \ , \  F = (W-T)Y \ , \ D = W-T,$$
il vient
$$EQ + FR  = D, \ \text{une identit\Ž de B\Žzout}$$
o\ 
$$\deg(Q) = q \ , \ \deg(R) = r \ , \  \deg(D)\leqs n-1 = q + r -1.$$
Dans ces conditions, on peut remplacer, dans l'identit\Ž de B\Žzout, les deux polyn\™mes, $E$ et  $F$, par des deux polyn\™mes, $G$ et $H$, respectivement,
avec $\deg(E) \leqs r-1$ et $\deg(F) \leqs q-1$. 

\

\lopar Le proc\Žd\Ž est classique. On effectue la division euclidienne de  $E$ par le polyn\™me  $R$ sous la forme
$E = CR + G \ , \ \deg(G) < \deg(R)$,
et on pose $H = C + F$. Il vient
$$(CR + G)Q + FR = GQ + (C+F)R = GQ + HR$$
autrement dit,
$GQ + HR = D$
avec $\deg(H) < \deg Q$.\}

\

\noi En poursuivant ainsi, on arrive \ˆ l'\Žgalit\Ž d\Žsir\Že, $UV = f(x)$. Dans $U$ et $V$, on regroupe les termes en $x^i$, pour $i = 0,1, \dots, n$. On reconstitue ainsi deux polyn\™mes, $g(x)$ et $h(x)$, pour lesquels on a $g(x)h(x) = f(x)$ et dont les coefficients sont de la forme (r\Žel + infinit\Žsimal), comme annonc\Ž !\qed

\

\head{Fin de la d\Žmonstartion du th\Žor\me}

\

Soit 
$$u(x) = x^n + a_1x^{n-1} + a_2x^{n-2} + \dots + a_{n-1} x + a_n,$$ 
un polyn\™me irr\Žductible, \ˆ coefficients dans le corps liminaire $K_\alpha$.  Une premi\re substitution, $x = y- a_1/n$,  ram\ne \ˆ un  polyn\™me de la forme 
$$v(y) =  y^n + a_2y^{n-2} + \dots + a_{n-1} y + a_n\ \text{o\} \ a_1 = 0.$$
Ce polyn\™me est irr\Žductible, \Žgalement. 
Deux cas seulement, pourraient de pr\Žsenter. 
\su{Premier cas}  Ce polyn\™me se r\Žduit au mon\™me $y^n$. Puisqu'il est irr\Žductible, c'est le polyn\™me $v(y) = y$, et tout est dit : $-a_1$ est alors une racine du polyn\™me $u(x)$.
\su{Second cas} Le polyn\™me $v(y)$ n'est pas r\Žduit \ˆ $y^n$.
La forme normale du coefficient $a_i$, pour $i\geqs 2$,  commence par un $\omega^{c_i}r_i$. On introduit le nombre $c  = \max\{c_i/i = 2\leqs i \leqs n\}$.  La substitution $y = \omega^c x$ conduit \ˆ un polyn\™me de la forme
$$w(x) = x^n + \sum_{i=2}^n \omega^{-ic}a_ix^{n-i}.$$
Par construction, on a $c_i/i \leqs c$, pour $i = 2, \dots, n$, avec \Žgalit\Ž pour l'un au moins des indices $i$. Le coefficient $\omega^{-ic} a_i$ de $x^{n-i}$ commence par $\omega^{-ic}\omega^{c_i} = \omega^{(c_i -ic)}$ o\ $c_i - ic\leqs 0$ est n\Žgatif mais \Žgal \ˆ $0$ pour un des indices au moins. Autrement dit, pour un $i\geqs 2$ au moins, la partie r\Želle du coefficicent $\omega^{-ic} a_i$,  a une partie r\Želle non nulle.

\

\noi Le polyn\™me $w(x)$ est irr\Žductible, comme le polyn\™me $u(x)$.
Il est de degr\Ž impair  et de la forme
$$w(x) = x^n  + \sum_{i = 2}^n (r_i + t_i)x^{n-i}$$
o\, pour un $i\geqs 2$ au moins, la partie r\Želle $r_i$ n'est pas nulle.
Il rel\ve du lemme de Hensel. La partie r\Želle, $P(x)$, du polyn\™me $w(x)$  ne se d\Žcompose donc pas en deux polyn\™mes premiers entre eux ! Ce polyn\™me, $P(x)$, est donc de la forme $(x-a)^n$ ou $(x^2 +bx +c)^m$. Mais la seconde forme est exclue puisque le polyn\™me est de degr\Ž impair.
De plus, le coefficient de $x^{n-1}$ dans ce polyn\™me est nul, donc $a=0$. Autrement dit, la partie r\Želle du polyn\™me est $x^n$ ce qui contredit le fait que l'un des $r_i$ n'est pas nul. Cela ach\ve la d\Žmonstration.\qed

\

\

{\bf Insistons. La classe des nombres, $\sS$, est un corps r\Žellement clos, r\Žunion de la famille croissante, $\{K_\alpha : \alpha \in \sO\}$, des corps liminaires, r\Žellement clos !}

\

\head{Encore une propri\Žt\Ž remarquable}

\

{\sc Conway} a d\Žmontr\Ž que tout corps totalement ordonn\Ž, K, o\ $K$ est un ensemble, poss\de une copie comme sous-corps du corps des nombres, $\sS$.  En ce sens, on peut dire que $\sS$ est un corps totalement ordonn\Ž {\bf universel} ! \emph{Mais cela est une autre histoire}, comme disait Kipling

\

 \head{Compl\Žments}
 
 \

  \head{C1 Les formes monomiales}
  
  \
  
  On a d\Žj\ˆ introduit les formes normales, les formes courtes et les formes longues des nombres. On va introduire une  quatri\me esp\ce de formes, les formes monomiales ! On appellera {\bf forme monomiale} toute famille de mon\™mes, $\{\omega^t u_t\}$, qui produit une forme normale lorsqu'on lui fait subir les op\Žrations d\Žcrites, ci-dessus, dans la d\Žfinition de la multiplication !
  
 \
 
 \noi Autrement dit,  une forme monomiale est une famille $\{\omega^t u_t\}$ de mon\™mes qui v\Žrifie les conditions suivantes. 
 
 \noi 1. Les $t$ sont des nombres, les $u_t$ des nombres r\Žels.

 \noi 2. L'ensemble des nombres $t$ est bien ordonn\Ž pour la relation $\geqs$.
 
  \noi 3. Pour chaque nombre $\tau$, il n'y a qu'un nombre fini de $t$ de la famille pour lesquels on a $t =\tau$.
 
 \
 
 \noi On dispose ainsi de {\bf quatre} esp\ces de formes pour \Žcrire les nombres. Cela pr\Žsente une souplesse tr\s estimable dans leur maniement !
 
 \
 
  \head{C2 Sur diff\Žrentes fa\c cons d'\Žcrire la multiplication}

 \
 
 On se donne deux nombres, $a$ et $b$, sous des formes longues :
 $$a = \sum_{i<\alpha} \omega^{a_i} r_i  \ , \ b = \sum_{j <\beta} \omega^{b_j}s_j.$$
 On pourra \Žcrire le produit sous forme longue :
  $$\left (\sum_{i <\alpha} \omega^{a_i}r_i\right).\left(\sum_{j < \beta} \omega^{b_j} s_j\right) = \left(\ \sum_{k < \gamma} \omega^{c_k} t_k\right).$$
On pourra \Žcrire aussi, sous {\bf forme monomiale} :
  $$\left (\sum_{i <\alpha} \omega^{a_i}r_i\right).\left(\sum_{j < \beta} \omega^{b_j}s_j\right) = \{\omega^{a_i + b_j}r_i. s_j\}.$$
 Enfin, plus simplement encore :
 $$\{\omega^{a_i}r_i\}. \{\omega^{b_j}s_j\}= \{\omega^{a_i + b_j}r_i. s_j\}.$$
 
 \

\head{C3 Utilisation des formes monomiales}

\

On peut multiplier, additionner, calculer des polyn\™mes et plus g\Žn\Žralement, des s\Žries enti\res de formes monomiales :

$x.y$

$ x + y$

$t_0 + t_1x + \dots + t_nx^n$ 

$t_0 + t_1x + \dots + t_nx^n + \dots = \sum_{n < \omega} t_nx^n$.

\su{Explications} On a d\Žj\ˆ vu comment on multiplie les formes monomiales. En particulier, cela permet de calculer les formes monomiales des puissances $x^n$ d'une forme monomiale $x$.

\

\noi Pour l'addition, la somme de deux formes monomiales est la r\Žunion disjointe de ces  deux formes :
$$\{\omega^{y_\alpha}.r_\alpha\}_{\alpha \in A}  + \{\omega^{z_\beta}.s_\beta\}_{\beta \in B} = \{\omega^{y_\alpha}.r_\alpha\} \uplus \{\omega^{z_\beta}.s_\beta\}.$$
\lopar Plus pr\Žcis\Žment, soit $C = A\uplus B$ la somme disjointe des deux ensembles d'indices, $A$ et $B$. Autrement dit, en rempla\c cant $A$ et $B$ par deux ensembles disjoints, $A'$ isomorphe \ˆ $A$, et $B'$ isomorphe \ˆ $B$, on prend $C = A' \cup B'$. On r\Žindexe les deux formes  monomiales \ˆ l'aide des indices $\gamma\in C$ et l'on prend
$$\{\omega^{y_\alpha}.r_\alpha\} \uplus \{\omega^{z_\beta}.s_\beta\} = \{\omega^{t_\gamma}.u_\gamma\}$$
ayant v\Žrifi\Ž, simplement, que c'est bien une forme monomiale !\}

\

\noi En particulier cela permet de calculer la forme monomiale d'une somme finie de formes monomiales, $x + y + \dots +z$. On saura aussi calculer la forme monomiale d'un polyn\™me $t_0 + t_1x + \dots + t_nx^n$ o\ $x$ est une forme monomiale.

\

\head{C4 Voici quelques propi\Žt\Žs compl\Žmentaires du corps des nombres}

\ 

Afin d'abr\Žger, on dira que le nombre $a = \sum \omega^{y_\alpha}r_\alpha$, sous forme longue, est {\bf tr\s positif} si tous les coefficients $r_\alpha$ sont positifs ! On d\Žsigne par $J$ la classe des nombres tr\s positifs. On introduit la relation binaire $a \succ b$, qui  dit $a - b$ est tr\s positif :
$$a \succ b \iff (a-b) \in J.$$
La relation $\succ$ est plus fine que l'ordre total $\geqs$.
$$a \succ b \implies a\geqs b.$$
Tout nombre tr\s positif est, en particulier, positif.

\

 \noi On v\Žrifie, simplement, que $J$ poss\de les propri\Žt\Žs suivantes.

\

\noi 1. On a  $J + J \inc J \ , \ J.J \inc J \ , \ \lambda. J +\mu.J\inc J$ pour tous nombres r\Žel positifs, $\lambda$ et $\mu$.

\noi C'est \ˆ dire que la classe $J$ est stable pour l'addition, la multiplication et le produit par les nombres r\Žels positifs. [C'est {\bf un c\™ne convexe} !]

\

\noi 2. On a $a \in J \cap (-J) \iff a = 0$.

\noi C'est \ˆ dire que le seul nombre $a\in J$ tel que $-a\in J$ est $a=0$.

\
 
 \noi 3. Pour chaque nombre $a = \sum \omega^{y_\alpha} r_\alpha$,
 on pose
$$a^+ = \sum_{r_\alpha \geqs 0}  \omega^{y_\alpha} r_\alpha \ , \ a^- = -\sum_{r_\alpha \leqs 0}  \omega^{y_\alpha} r_\alpha.$$
Ce sont deux nombres  tr\s positifs. On a
 $$a = a^+ -a^- \ , \ a^+\in J \ , \  a^-\in J.$$
 On a ainsi  $J-J = \sS$. Tout nombre est diff\Žrence de deux nombres tr\s positifs.  
 
 \
 
\noi La relation $\succ$ est une relation d'ordre sur la classe des nombres. En effet, l'on a 

$a-b \in J \et b-c\in J \implies a-c \in J$

$a-b \in J  \et b-a\in J \implies a=b$

ce qui veut dire

$a \succ b \et  b\succ c \implies a\succ c$

$a\succ b \et b\succ a \implies a=b$. 

\noi Autrement dit cette relation est transitive, antisym\Žtrique et r\Žlexive. C'est bien une relation d'ordre. 

\

\noi  Pour abr\Žger, on dira que $a$  est {\bf bien plus grand} $b$,  et $b$ est {\bf bien plus petit} que $a$, lorsque l'on a $a\succ b$.

\

\noi Cet ordre n'est pas un ordre total. Il y a des couples de nombres incomparables ! Ainsi, par exemple, $(\omega -1)$ n'est ni bien plus grand ni bien plus petit que $0$ ! Cependant, tout couple, $a,b$, poss\de une borne sup\Žrieure, $a\lor b$, et une borne inf\Žrieure, $a\land b$, \ˆ savoir, 
pour
$$a = \sum\omega^{y_\alpha} r_\alpha \ , \ b = \sum\omega^{y_\alpha} s_\alpha,$$
 $$a\lor b = \sum \omega^{y_\alpha}({r_\alpha} \lor {s_\alpha}) \ , \ a\land b = \sum\omega^{y_\alpha}({r_\alpha} \land {s_\alpha}).$$
 
\
 
 \noi La d\Žcomposition de $a$ comme diff\Žrence de deux nombres tr\s positifs est bien entendu unique. On pose
 $$|a| = a^+ + a^-.$$
 On aura $$|a|=  a\lor (-a).$$
 C'est un nombre tr\s positif que l'on appelle {\bf la valeur absolue} du nombre $a$. Par exemple, on a $0\lor(\omega -1)= (\omega -1)\lor (1-\omega) = \omega + 1$. Les deux nombres $0$ et  ($\omega - 1$) ne sont pas comparables !
 
 \
 
 \noi On a aussi 
  $$|a|=  a\lor (-a) \ , \ |-a| = |a| \ , \ |a| - |b| \prec |a+b| \prec |a| + |b|,$$
  On notera encore cette double in\Žgalit\Ž qui en d\Žcoule :
  $$|a| - |b| \leqs |a+b| \leqs |a| + |b|.$$
 Tout cela est bien connu pour les groupes ordonn\Žs {\bf r\Žticul\Žs} !
 
 \

\

\head{\'Epilogue}

\

On a montr\Ž comment pr\Žsenter simplement la classe, $\sS$,  des nombres de {\sc Cuesta-Conway} en partant de leurs formes normales qui sont des sommes de familles d'exponentielles ind\Žx\Žes par des ordinaux. On a montr\Ž comment on y d\Žfinit addition, multiplication et relation d'ordre totale et que $\sS$ est alors un corps commutatif totalement ordonn\Ž et r\Žellement clos !

\

\noi Les citations suivantes de {\sc Gonshor} montrent comment ce dernier avait pressenti la possibilit\Ž d'introduite les nombres surr\Žels {\it \ˆ l'envers}, \ˆ notre mani\re, en se passant des complexit\Žs de la construction donn\Že par {\sc Conway}. Mais il ne l'a pas fait !

\ 

\noi {\sc Gonshor} [ 4 ] p. 71, \Žcrit : {\sl \lq\lq Theorems 5.7 and 5.8 give us a powerful tool for dealing with surreal numbers. In fact, for many purposes we can simply work with these generalized power series and ignore what surreal numbers are in the first place. This is an example of the whole spirit of abstraction in mathematics.'' }

\

\noi Les th\Žor\mes 5.7 et 5.8 traitent de l'addition et de la multiplication de nombres \ˆ l'aide des formes longues et monomiales !

\

\noi Il \Žcrit p.71, plus bas, \lq\lq {\sl There is another method of using generalized power series to obtain existence results which does not depend on familiarity with identities for ordinary formal series ... I like this method because of its elementary self-contained algebraic nature.''}

\

\noi C'est ce que {\sc Harry  Gonshor}, [4, p.71], \Žcrit en parlant des sommes formelles index\Žes par les ordinaux. 

\

\noi Plus bas encore, page 113, il dit ceci : \lq\lq 
{\sl In fact, many of the arguments depend only on the generalized group ring structure given by the normal form and make no essential use of the fact that the exponents are surreal numbers. Thus by taking theorems 5.7 and 5.8 as definitions for a set consisting of generalized power series we can obtain a theory which is independent of the theory of surreal numbers although the results will be valid for surreal numbers as one special case. One convenient hypothesis for the exponents is that they form a divisible abelian group.''}

\

\noi Pourtant, je ne connais aucun auteur qui ait eu d\Žj\ˆ l'id\Že de proc\Žder de la m\me mani\re en partant des formes normales. Si, cependant, il y en a qui l'on fait, ils pourraient me le signaler. D'avance, je leur pr\Žsente mes excuses de ne pas les avoir cit\Žs, ne les connaissant pas !

 \
 
 \head{Petit lexique}
 
 \ 
 
 doublure 7

 coupure de Dedekind 1
 
 corps liminaire 10
 
 exponentielle 2

forme normale 4

forme courte 4

formes longues 4

forme monomiale 15

hauteur 3

infinit\Žsimal 12

nombre 1 

nombres de Cuesta-Conway 1

nombres surr\Žels 1

partie r\Želle 12

\

\head{Biliographie}

\

\noi[ 1 ] {\sc N. L. ALLING}, {Foundations of analysis over surreal number fields}, Mathematics Studies 141, xvi + 373 pp., North Holland, 1987.

\

\noi[ 2 ] {\sc J. H. CONWAY}, {On numbers and games}, ix + 230 pp.,
Academic Press Inc. 1976, reprinted 1979. 

\

\noi[ 3 ] {\sc N. CUESTA}, {\sl Algebra ordinal}, Revista de la Real Academia De Ciencias Exactas, Fisicas Y Naturales, {\bf 58} no 2 (1954) 103-145.

\

\noi[ 4 ] {\sc H. GONSHOR}, {\sl An introduction to the theory of surreal numbers}, London Mathematical Society Lecture Note Series. 110, Cambridge University Press, 1987, Re-issued in digitally printed version 2008, ii + 192 pages.

\

\noi[ 5 ] {\sc L. HADDAD}, {\sl  A propos des nombres de Conway : lettre \ˆ un ami}, {\tt arXiv:2101.11790v2}, 30 Jan 2021.

\

\noi[ 6 ] {\sc L. HADDAD}, {\sl Un lemme combinatoire de H. B. Neumann}, {\tt	arXiv:2203.01119}, 1 Mar 2022. 

\

\noi[ 7 ] {\sc S. MacLANE}, {\sl The universality of formal power series fields},  Bull. A. M. S., {\bf 45} (1939) 888-890.

\

\noi[ 8  ] {\sc B. H.  NEUMANN, {\sl On ordered division rings.} Trans. Amer. Math. Soc. {\bf 66} (1949) 202-252.

\

\

\

\

\enddocument